\documentclass{amsart}\usepackage{times}
\usepackage{amssymb,amscd}
\usepackage[11pt,curve,matrix,arrow,frame,graph]{xy}
\newtheorem{thm}{Theorem}[section]
\newtheorem{lem}[thm]{Lemma}
\newtheorem{prop}[thm]{Proposition}
\newtheorem{cor}[thm]{Corollary}

\theoremstyle{remark}
\newtheorem{rem}[thm]{Remark}
\newtheorem*{merci}{Acknowledgements} 

\numberwithin{equation}{section}   

\newcounter{counteroman}
\newenvironment{enumeroman}{\begin{list}{\roman{counteroman})}{\usecounter{counteroman}}}{\end{list}}

\renewcommand{\H}{\mathbb{H}}
\newcommand{\R}{\mathbb{R}}
\newcommand{\bS}{\mathbb{S}}\newcommand{\N}{\mathbb{N}}\newcommand{\Z}{\mathbb{Z}}\newcommand{\V}{\mathbb{V}}

\newcommand{\mcH}{\mathcal{H}}

\newcommand{\mcO}{\mathcal{O}}\newcommand{\mcV}{\mathcal{V}}

\newcommand{\la}{\langle}\newcommand{\ra}{\rangle}

\newcommand{\bpm}{\begin{pmatrix}}
\newcommand{\epm}{\end{pmatrix}}

\DeclareMathOperator{\im}{Im}

\DeclareMathOperator{\rank}{rank}\DeclareMathOperator{\Card}{Card}

\DeclareMathOperator{\support}{support}
\DeclareMathOperator{\tr}{tr}
\DeclareMathOperator{\vol}{vol}

\title{Rigidity and $L^2$ cohomology of hyperbolic manifolds.}
\author{G. Carron}

\begin{document}

\maketitle
\begin{minipage}{10cm}{\footnotesize ABSTRACT: 
When $X=\Gamma\backslash \H^n$ is a real hyperbolic manifold, it is already known that  if  the critical exponent is 
 small enough then some cohomology spaces and some spaces of $L^2$ harmonic forms vanish. In this paper, we show rigidity 
results in the borderline case of these vanishing results. 
}\end{minipage}
\vskip0.2cm
\begin{minipage}{10cm}
{\footnotesize R\'ESUM\'E : La petitesse de l'exposant critique du groupe fondamental d'une vari\'et\'e hyperbolique 
implique des r\'esultats  d'annulation pour certains espaces de cohomologie et de formes harmoniques $L^2$.
 Nous obtenons ici des r\'esultats de rigidit\'e reli\'es \`a ces r\'esultats d'annulations.
  Ceci est une g\'en\'eralisation de r\'esultats d\'ej\`a connus dans le cas convexe co-compact.
}
\end{minipage}

\section{Introduction.}
When $\Gamma$ is a discrete torsion free subgroup of $\mathrm{SO}(n,1)$,  its critical exponent is defined by 
$$\delta(\Gamma):=\inf\{s>0, \sum_{\gamma\in \Gamma}e^{-sd(z,\gamma(z))}<+\infty\}.$$
It is easy to see that this definition doesn't depend on the choice of $z\in \H^n$ and that 
$\delta(\Gamma)\le n-1$. The critical exponent measures the growth of $\Gamma$-orbits:
$$\delta(\Gamma)=\limsup_{R\to +\infty} \frac{\log\Card (B(z,R)\cap \Gamma.z)}{R}.$$
An important and beautiful result of D. Sullivan \cite{sullivanl} (see also \cite{patterson}  in dimension $2$)
 is a formula between the critical exponent  and the bottom of  the spectrum of the Laplacian (on functions) on the manifold
  $\Gamma\backslash \H^n$ :

\begin{center}
\noindent\textsl{If $\delta(\Gamma)\le (n-1)/2\ $ then $\lambda_0\left(\Gamma\backslash \H^n\right)=(n-1)^2/4$.}

\noindent  \textsl{If $\delta(\Gamma)\ge (n-1)/2\ $ then $\lambda_0\left(\Gamma\backslash \H^n\right)=\delta(n-1-\delta)$.}
\end{center}
When $\Gamma$ is geometrically finite, the critical exponent is linked with the Hausdorff dimension of the limit 
set $\Lambda(\Gamma)=\overline{\Gamma.z}\cap \partial_\infty \H^n$ (where the closure is taken in the geodesic
 compactification of the hyperbolic space
$  \overline{ \H^n}=\H^n\cup   \partial_\infty \H^n$) or with the entropy of the geodesic flow
(\cite{sullivand},\cite{patterson},\cite{bishop}).

When $\Gamma$ is convex-cocompact and isomorphic to a cocompact discrete torsion free subgroup of $\mathrm{SO}(n-p,1)$, then
$$\delta(\Gamma)\ge n-1-p$$ with equality if and only if $\Gamma$ stabilizes cocompactly a
 totally geodesic $(n-p)$-hyperbolic subspace $\H^{n-p}\subset \H^n$ (\cite{bowen},\cite{bourdon},\cite{yue},\cite{BCGacta}).
  Other rigidity results in terms of the critical exponent have been recently obtained for amalgamated products (\cite{shalom},\cite{BCGcras},\cite{BCGjdg}).

Using different techniques, H. Izeki, H. Izeki and S. Natayani (\cite{izeki},\cite{IzekiN}) and X. Wang \cite{wang2} obtained 
rigidity results based on the De Rham cohomology with compact support
 \footnote{or cohomology in complementary degree using Poincar\'e's duality} :
\begin{thm}\label{cc} Let $X=\Gamma\backslash \H^n$ be  a convex-cocompact hyperbolic manifold, assume that for some $p<n/2$ :
$$H^p_c(X)\not=\{0\}$$
then 
$$\delta(\Gamma)\ge n-1-p$$ 
with equality if and only if 
 $\Gamma$ stabilizes cocompactly  and positively a totally geodesic $(n-p)$-hyperbolic subspace $\H^{n-p}\subset \H^n$.
\end{thm}
We say that $\Gamma$ \textsl{stabilizes cocompactly  and positively} a totally geodesic $(n-p)$-hyperbolic 
subspace $\H^{n-p}\subset \H^n$ when  $\Gamma$ stabilizes cocompactly a totally geodesic $(n-p)$-hyperbolic subspace 
$\H^{n-p}\subset \H^n$ and when $\Gamma$ acts trivially on the orientation normal bundle of $\H^{n-p}\subset \H^n$ .

In fact, R. Mazzeo has shown that  the cohomology with compact support of a convex cocompact hyperbolic $n-$manifold
 is isomorphic to the space of harmonic $L^2$ forms in degree $p<n/2$ \cite{Mazzeo}: 
If $X=\Gamma\backslash \H^n$ is convex cocompact and if $p<n/2$ then
$$H^p_c(X)\simeq \mcH^p(X):=\{\alpha\in L^2(\Lambda^pT^*X), d\alpha=d^*\alpha=0\}.$$

In \cite{CP}, with E. Pedon, we obtained the following result : 
\begin{thm}Let $X=\Gamma\backslash \H^n$ be  a hyperbolic manifold, assume that for $p<n/2$ :
$$\delta(\Gamma)< n-1-p,$$  then 
$$\mcH^p(X)=\{0\}.$$ Moreover the bottom of the spectrum of the Hodge-De Rham Laplacian on $p$ forms  is bounded from below by :
$$\lambda_0(dd^*+d^*d,\Gamma\backslash \H^n)\ge (\delta(\Gamma)-p)(n-1-p-\delta(\Gamma)) $$ 
if $\ (n-1)/2\le \delta(\Gamma)<n-1-p$ and  
 $$\lambda_0(dd^*+d^*d,\Gamma\backslash \H^n)\ge \frac{(n-1-p)^2}{4}$$
  if $\delta(\Gamma)\le (n-1)/2.$
\end{thm}
Together with Mazzeo's interpretation of the space of $L^2$ harmonics forms, this result implies a part of the theorem 
\ref{cc} : convex cocompact hyperbolic manifolds with non trivial cohomology with compact support in some degree $p<n/2$
 have a critical exponent strictly larger than $n-1-p$. 
In this paper, we study rigidity result without the convex cocompact hypothesis. 
The case of degree $p=1$ is covered by the following very general result of P. Li and J. Wang \cite{LiWang} 
(see also \cite{wang1} for earlier results):
\begin{thm}\label{LiW}
 If $(M^{n>2},g)$ is a complete Riemannian manifold with
${\rm Ricci\,}_g\ge -(n-1)g$ and
$\lambda_0(M^n,g)\ge (n-2)$, then either
\begin{enumeroman}
\item  $M$ has only one end with infinite volume or 
\item $(M^n,g)$ is isometric to the warped product $$(\R\times N, (dt)^2+\cosh^2(t) h)$$ with 
$(N,h)$ compact and ${\rm Ricci\,}_h\ge -(n-2)h$  .
\end{enumeroman}\end{thm}
This has the following consequence for hyperbolic manifold :
\begin{thm} Let $X=\Gamma\backslash \H^n$ be a hyperbolic manifold with $n>2$ and whose injectivity radius  is positive.
If 
$$\delta(\Gamma)\le n-2$$ then either
\begin{enumeroman}
\item  $H^1_c(X,\Z)=\{0\}$ or
\item $H^1_c(X,\Z)\not=\{0\}$, $\delta(\Gamma)= n-2$ and  $\Gamma$ stabilizes cocompactly 
a totally geodesic hypersurface $\H^{n-1}\subset \H^n$
\end{enumeroman}
\end{thm}

Recently M. Kapovich has studied the link between the critical exponent and the cohomological dimension of 
a hyperbolic manifold  $X=\Gamma\backslash \H^n$ relative to the $\epsilon$-ends
whose rank is larger or equal to $2$ \cite{kapovich}. Let $\epsilon$ be a positive number
 smaller that the Margulis constant and denote $X_{<\epsilon}$  the $\epsilon$-thin part of $X$,
  that is the set of point of $X$ where the injectivity radius is smaller than $\epsilon$. 
  Let $M_{<\epsilon}\subset X_{<\epsilon}$ be the union of the connected component of $X_{<\epsilon}$ whose fundamental 
  group has rank larger than or equal to $2$. M. Kapovich proves the following result :
\begin{thm} Assume that $$\delta(\Gamma)< n-p-1$$
and let $R$ be a commutative ring with unit and $\V$ be a $R\Gamma$-module
then  $$H^{n-p}(X, M_{<\epsilon},\V)=\{0\}.$$
Moreover, when $\Gamma$ is  assumed to be geometrically finite and when 
 $$\delta(\Gamma)= n-p-1\ \mathrm{and}\  H^{n-p}(X, M_{<\epsilon},\V)\not=\{0\},$$ then 
$\Gamma$ stabilizes  a totally geodesic $n-p$-hyperbolic subspace $\H^{n-p}\subset \H^n$ with 
$\vol (\Gamma\backslash \H^{n-p})<\infty$.
\end{thm}

We'll prove a similar rigidity result for the De Rham cohomology in the space of $L^2$ harmonic form :

\begin{thm} Let $X=\Gamma\backslash \H^n$ be a oriented hyperbolic manifold with $n>3$ and
 let $E\subset X_{<\epsilon}$ be the union of all unbounded connected components of the $\epsilon$-thin part and let $p<(n-1)/2$.
There is a natural linear map 
$$\mcH^{n-p}(X)\rightarrow H^{n-p}(X,E).$$
Moreover if
$\delta(\Gamma)=n-p-1$ and if this map is non zero then $\Gamma$ stabilizes positively a totally geodesic 
$n-p$-hyperbolic subspace $\H^{n-p}\subset \H^n$ with 
$\vol (\Gamma\backslash \H^{n-p})<\infty$.
\end{thm}
We recall that when $\delta(\Gamma)<n-p-1$ and $p<n/2$, then $\mcH^{n-p}(X)=\{0\}$.
The above restriction on the range of the degree $p<(n-1)/2$ comes from the fact that when $p\ge (n-1)/2$, 
we are not able to build a map $\mcH^{n-p}(X)\rightarrow H^{n-p}(X,E)$, however we'll give a similar result based on the cohomology
with compact support  (see theorem \ref{compact}).

When $\Gamma$ is geometrically finite,
 then a quick look at the topological interpretation of the space of $L^2$ harmonic forms obtained
  by R. Mazzeo and R.Phillips \cite{MP} shows that 
we have in this case :
$\mcH^{n-p}(X)\simeq H^{n-p}(X,E)$ (see \S \ref{GF}).
Moreover using the proof of \cite[lemma 8.1]{kapovich}, when $n-p>1$
we obtained 
$H^{n-p}(X,E)=H^{n-p}(X,M_{<\epsilon})$. Hence in the geometrically finite case, 
we are able to recover the rigidity result of M. Kapovich for $R=\R=\V$.

We now describe the proof of our result. Our proof owns a lot to X. Wang's proof of the theorem \ref{cc}
 but with several new technical points.
 
The first point is  to analyse the case of equality in the refined Kato inequality.
When $p<n/2$ and when $\xi$ is a harmonic $L^2$ $p$-form on the hyperbolic manifold $X=\Gamma\backslash \H^n$,
 then the refined  Kato's inequality (\cite{branson}\cite{kato}) implies that the function 
$$\phi:=|\xi|^{\frac{n-1-p}{n-p}}$$
satisfies
\begin{equation}\label{subeigen}\Delta \phi\le p(n-1-p) \phi.
\end{equation}
Our first result describes the equality case in this inequality; 
it is an extension of a result of X. Wang who described the equality case when
 $\xi$ is squared integrable and $\Gamma$ is convex cocompact.

Then we describe how we can define a map
$$\mcH^{n-p}(X)\rightarrow H^{n-p}(X,E)$$ or
$$H^p_c(X\setminus E)\rightarrow \mcH^p(X)$$
when $p<(n-1)/2$.  Note that $E$ being an open set, $X\setminus E$ is a closed subset of $X$ 
and that forms with compact support in $X\setminus E$ have a support that can touch $\partial E$.

The second crucial point is to prove that when the map $$H^p_c(X\setminus E)\rightarrow \mcH^p(X)$$ is not zero and 
$\delta(\Gamma)=n-1-p$ then there is a $L^2$ harmonic $p$-form  $\xi$ such that 
$\phi:=|\xi|^{\frac{n-1-p}{n-p}} \in L ^2$.

Then according to D. Sullivan's result, the bottom of the spectrum of the Laplacian on $X$ is 
$\delta(\Gamma) (n-1-\delta(\Gamma))=p(n-1-p)$, it is easy to deduce that 
in fact $\phi$ is a eigenfunction of the Laplace operator. Then we use our description of the equality in (\ref{subeigen}).

In the degree $p=(n\pm 1)/2$, then our methods does not apply because we are not able to build a map
$H^p_c(X\setminus E)\rightarrow \mcH^p(X)$.
However, there is always a map 
$H^p_c(X)\rightarrow \mcH^p(X)$
and our  proof will also show the following result :
\begin{thm}\label{compact}Let $X=\Gamma\backslash \H^n$ be a hyperbolic manifold with $n>3$.
Assume that  for a $p<n/2$:
$$\delta(\Gamma)\le n-1-p\ ,$$
then the image of the cohomology with compact support in the absolute cohomology is zero in degree $p$ :
$$\im(H_c^p(X)\rightarrow H^p(X))=\{0\}.$$
Moreover either 
\begin{enumeroman}
\item The map $H^p_c(X)\rightarrow \mcH^p(X)$ is zero 
\item Or the map $H^p_c(X)\rightarrow \mcH^p(X)$ is an isomorphism, $\delta(\Gamma)=n-1-p$ and
 $\Gamma$ stabilizes cocompactly and positively a totally geodesic $(n-p)$-hyperbolic subspace $\H^{n-p}\subset \H^n$.
\end{enumeroman}
\end{thm}

\begin{rem}
The case of hyperbolic manifolds of dimension $3$ is already covered by P. Li and J. Wang's result \ref{LiW}.
\end{rem}

\begin{merci}This text is an attempt to answer some of the  questions that have been asked after my talk at the 
conference "Spectral Theory and Geometry"  in honour of my teachers, advisor : P. B\'erard and S. Gallot.
 Hence it is a pleasure to thank G. Besson, L. Bessi\`eres, Z. Djadli for having organized this very nice conference. 
  I also take the opportunity to thank P. B\'erard and S. Gallot for all the beautiful mathematics that they taught  me. 
I thank V. Minerbe for his useful comments. I was partially supported by the project ANR project GeomEinstein 06-BLAN-0154.
  Eventually, I want to dedicate my paper to H. Pesce.
\end{merci}

\section{The equality case in the refined Kato's inequality on hyperbolic space.}
The classical Kato inequality says that if $\xi$ is a smooth $p$-form on a Riemannian manifold
$(M^n,g)$, then 
$$\left|d|\xi|\right|^2\le \left| \nabla \xi\right|^2.$$
When $\xi$ is assumed to be moreover closed and coclosed :
$$d\xi=d^*\xi=0,$$
then this Kato inequality can be refined :
\begin{equation}
\label{eq:1}
\frac{n+1-p}{n-p}\left|d|\xi|\right|^2\le \left| \nabla \xi\right|^2 ,
\end{equation}
See \cite{bourguignon}, for a convincing explanation of the principle leading to this inequality, and \cite{branson},
 \cite{kato} for the 
computation of the refined Kato constant. Where $(M^n,g)$
is a hyperbolic manifold, direct computations show that if $\xi$ is a harmonic
$p$-form then the function $\phi:=|\xi|^{\frac{n-1-p}{n-p}}$ satisfies
$$\Delta\phi\le p(n-1-p)\phi,$$
with equality if and only if we have equality in the refined Kato inequality (\ref{eq:1}).
In this situation, X. Wang has described the equality case when $\xi$ is square integrable 
and when $(M,g)$ is convex cocompact \cite{wang2}. 
Our first preliminary result is an extension of this result of X. Wang :

\begin{thm}\label{thm:Kato} Let $p,n$ be integers with $2p<n$ and $n>3$.
 If $\xi$ be a harmonic $p-$form on $\H^n$  such that we have everywhere equality :
$$\left| \nabla \xi\right|^2=\frac{n+1-p}{n-p}\left|d|\xi|\right|^2,$$
 then either
\begin{enumeroman}
\item there is a real constant $A$, an isometry $\gamma$ and a parallel $(p-1)-$ form $\omega$ on $\R^{n-1}$ such that
in the upper-half-space model of the hyperbolic space $$\H^n\simeq \{(y,x)\in (0,+\infty)\times\R^{n-1} \}$$ 
endowed with the Riemannian metric $\frac{(dx)^2+(dy)^2}{y^2}$, we have
$$\gamma^*\xi=A y^{n-1-p}\pi^*\omega,$$
where $\pi(y,x)=x$, or
\item there is a totally geodesic copy $\H^{n-p}\subset \H^n$ such that in Fermi coordinates around 
this $\H^{n-p}$ $$\H^n\setminus \H^{n-p}\simeq (0,+\infty)_t\times \bS^{p-1}\times \H^{n-p}$$
we have
$$\xi=A\frac{(\sinh t)^{p+1}}{(\cosh t)^{n-p+1}}dt\wedge d\sigma ;$$
recall that here $t$ is the geodesic distance to $\H^{n-p}\subset \H^n$ and $d\sigma$ is the volume form of $\bS^{p-1}$.
\end{enumeroman}
\end{thm}
Our arguments will follow closely those of X. Wang, however in his situation  only the case  ii) appears.

\subsection{Proof of the theorem \ref{thm:Kato}} Let $p,n$ be integers with $p<n/2$ and $n>3$. 
We consider $\xi$  a non trivial harmonic $p-$form on $\H^n$ such that we have everywhere: 
$$\left| \nabla \xi\right|^2=\frac{n+1-p}{n-p}\left|d|\xi|\right|^2.$$
Then $\phi:=|\xi|^{\frac{n-1-p}{n-p}}$ satisfies
$$\Delta\phi=p(n-1-p)\phi.$$
Then the Harnack inequalities imply that $\phi$ is positive.

According to \cite{kato}, there is locally a $1$-form $\alpha$ such that
$$\alpha\wedge \xi=0$$
and
$$\nabla \xi =\alpha\otimes \xi-\frac{1}{n+1-p}\sum_j \theta^j\otimes\theta^j\wedge \alpha^\sharp \llcorner\xi$$
for a local orthonormal dual frame $(\theta^1,...,\theta^n).$

Now we let $k:=n-1-p$ so $\phi=|\xi|^{\frac{k}{k+1}}.$
If $X$ is a vector field, then
\begin{equation*}
\begin{split}
\nabla_X\phi=&\frac{k}{k+1}|\xi|^{-\frac{1}{k+1}-1}\la\nabla_X\xi,\xi\ra\\
&=\frac{k}{k+1}|\xi|^{-\frac{1}{k+1}-1}\left[\alpha(X) |\xi|^2-\frac{1}{k+2}
\la  \alpha^\sharp \llcorner\xi,X\llcorner\xi\ra\right]
\end{split}
\end{equation*}
But we have $\alpha\wedge \xi=0$ hence
$$\la  \alpha^\sharp \llcorner\xi,X\llcorner\xi\ra=\alpha(X) |\xi|^2$$
and we obtain
\begin{equation}
\label{phi}
\nabla\phi=\frac{k}{k+2}\phi\, \alpha
\end{equation}
Hence, $\nabla \phi$ vanishes only where $\alpha$ vanishes and $\alpha$ is a smooth $(p-1)$ form.
We work on the open set 
$$U:=\{z\in \H^n, \nabla\phi(z)\not=0\}.$$
On $U$, we can locally find a orthonormal dual frame
$(\theta^1,...,\theta^n)$ such that
$$\alpha=(k+2)u \,\theta^1$$
with $u>0$.
Hence
\begin{equation}
\label{phi1}
\nabla\phi=ku \phi\, \theta^1
\end{equation}
As $\alpha\wedge\xi=0$, we can locally find a $(p-1)$-form $\omega$ such that
$$\xi=\theta^1\wedge \omega.$$
And we have
$$\nabla \xi=u\left[(k+1)\theta^1\otimes \theta^1\wedge\omega-\sum_{j=2}^n \theta^j\otimes\theta^j\wedge\omega \right]$$
Let $(e_1,...,e_n)$ be the frame dual to $(\theta^1,...,\theta^n)$, then we  obtain
$$\nabla_{e_1}\xi=(k+1)u\theta^1\wedge\omega=\nabla_{e_1}\theta^1\wedge\omega+\theta^1\wedge\nabla_{e_1}\omega$$
and for $j>1$ :
$$\nabla_{e_j}\xi=-u\theta^j\wedge\omega=\nabla_{e_j}\theta^1\wedge\omega+\theta^1\wedge\nabla_{e_j}\omega$$
With the fact that $\la \nabla_{e_j}e_1,e_1\ra=0$ and $\la \nabla_{e_j}\theta^1,\theta^1\ra=0$, 
we get the following identities
\begin{equation}
\label{11}
\theta^1\wedge\left[\nabla_{e_1}\omega-(k+1)u\omega\right]=0
\end{equation}
\begin{equation}
\label{12}
\nabla_{e_1}\theta^1\wedge\omega=0
\end{equation}
\begin{equation}
\label{21}
\theta^1\wedge\nabla_{e_j}\omega=0
\end{equation}
\begin{equation}
\label{22}
\left(\nabla_{e_j}\theta^1+u\theta^j\right)\wedge\omega=0
\end{equation}
Let $c\in \phi(U)$ and let $\Sigma_c:=U\cap \phi^{-1}\{c\}$, 
this is a smooth hypersurface and $e_1$ is an unit  normal vector field to $\Sigma_c$.
Then the equality (\ref{21}) implies that the pull back of $\omega$ to $\Sigma_c$ is parallel. 

At $z\in \Sigma_c$, we decompose
\begin{equation}
\label{decomp}
T_z\Sigma_c=E_z\oplus E_z^\perp
\end{equation}

where 
$$E_z:=\{v\in T_z\Sigma, v^\flat\wedge\omega=0\}$$
Let $L$ be the shape operator of $\Sigma_c$ at $z$ 
$$L\,:\,T_z\Sigma_c\rightarrow T_z\Sigma_c$$
$$Lv=-\nabla_ve_1;$$
we have
$$\nabla_{e_j}\theta^1=-\sum_{i=2}^n\la Le_j,e_i\ra\theta^i.$$
The equation (\ref{22}) implies that $L(E_z)\subset E_z$. Since
$L$ is a self adjoint operator we also have 
$L\left(E_z^\perp\right)\subset E_z^\perp$ and moreover still according to equation (\ref{22}), we have
$$L(X)=uX,\forall X\in E_z^\perp.$$
Then $\omega$ being parallel, the decomposition (\ref{decomp}) induced a parallel decomposition of the tangent bundle of $\Sigma_c$, in
particular if $X\in E_z, Y\in E_z^\perp$ are unit vectors
then the sectional curvature of $\Sigma_c$ in the direction of $X\wedge Y$ is zero and the Gauss Egregium theorem implies that
$$-1=\left(\la LX, Y \ra\right)^2-\la LX,X\ra \la LY,Y\ra$$
hence we have
$$LY=\frac{1}{u} Y,\  \forall Y\in E_z.$$ 
We can now compute the Ricci curvature of $\Sigma_c$, it is given by the formula
$${\rm ricci}_{\Sigma_c}=\left(\rank E_z-1\right)\left(-1+\frac{1}{u^2}\right)g_{E_z}+\left( \rank E_z^\perp-1\right) (-1+u^2)g_{E_z^\perp}$$
The hypothesis $n-1>2$ and the trace of the Bianchi identity
$$\delta_{g_{\Sigma_c}} {\rm ricci}_{\Sigma_c}=-d{\rm Scal}_{\Sigma_c}$$
implies that the function $u$ is constant 
on each connected component of $\Sigma_c$.
As $u$ is proportional to the length of $\nabla\phi$, this implies that
$$\forall j>1,\ 0=\la \nabla_{e_j} \nabla \phi, e_1\ra=\la e_j, \nabla_{e_1}\nabla \phi\ra . $$
So that $ \nabla_{e_1}e_1=0$, and $\theta^1$ is (locally and up to a sign) the differential
of the distance to $\Sigma_c$ and $\phi$ is a function of the sign distance to $\Sigma_c$.

\noindent{\it First case : $u=1$ at a point $z$:}  then the connected component of $\Sigma_c$ which contains $z$ is a totally umbilical flat hypersurface of $\H^n$.
Up to an isometry, we can assume that this connected component of $\Sigma_c$ 
is included in the horosphere :
$$\{y=1\}$$
in the upper half-space model of the hyperbolic space.

The facts that $\phi$ depends only on the distance to $\Sigma_c$ and that $\phi$ is a eigenfunction of the Laplace operator 
imply that there are constants $A, B$ such that in a neighborhood of $z$, we have
$$\phi(y,x)=Ay^{n-1-p}+By^{p}$$
In this case we have
$$\left|\nabla\phi\right|=\left|A(n-1-p)y^{n-1-p}+Bpy^{p}\right|.$$
But  $p\le (n-1)/2$, hence 
$$\left|\nabla\phi\right| \le (n-1-p)\phi$$
with equality at a point if and only if $B=0$.
But at $z$, we have $u=1$ hence (cf. \ref{phi1}) at $z$ we have
 $$\left|\nabla \phi\right|=(n-1-p)\phi.$$
 So that $B=0$ 
and $u=1$ around $z$. The unique continuation property for eigenfunctions of the Laplace operator implies that we have everywhere :
$$\phi(y,x)=Ay^{n-1-p}.$$ 
Hence we have also : $U=\H^n$, $e_1=y\partial_y$ and that $\omega$ is a parallel $(p-1)$-form on each horosphere
$\{y=c\}$. Finally, the equation
$$\nabla_{e_1}\omega=(k+1) \omega$$ implies that for a certain
 $$\tilde\omega\in \Lambda^{p-1}\left(\R^{n-1}\right)^*$$
we have
$$\xi =Ay^{n-1-p} dy\wedge \pi^*\tilde\omega$$
where $\pi(y,x)=x$.

\noindent{\it Second case : $u=1$ nowhere.}
The distributions induced by $E_z$ and $E_z^\perp$ are parallel hence integrable. Locally
there is a splitting
$$\Sigma_c=\Sigma_c(E)\times \Sigma_c\left(E^\perp\right).$$
And each $ \Sigma_c(E)$ has curvature $-1+u^{-2}$ and
 each $ \Sigma_c\left(E^\perp\right)$ has curvature $-1+u^{2}$ .
 We have $\omega=\Omega_E\wedge \tau$, where $\Omega_E$ is the volume form of $\Sigma_c(E)$
 and  $\tau$ is a parallel form on $\Sigma_c\left(E^\perp\right)$; however, the curvature of
 $\Sigma_c\left(E^\perp\right)$ is constant, not zero hence $\Sigma_c\left(E^\perp\right)$ has only
 parallel form in degree $0$ or in degree $\dim \Sigma_c\left(E^\perp\right)$.
 This implies that $\deg \omega=\rank E=p-1$.
 
We fixed now $z_0\in U$ and $c_0=\phi(z_0)$.
 We consider a neighborhood $\mcO$ of $z_0$ such that
 $\mcO\cap \Sigma_{c_0}$ is connected and isometric to 
 $S\times T$ where $S$ has curvature $-1+u^{-2}$ and $T$ has curvature $-1+u^{2}$.
 This neighborhood can be choosen so that the exponential map

 \begin{equation*}\begin{split}
E\,:\, (-\delta,\delta)\times \Sigma_c\cap \mcO &\rightarrow \mcO\\
(t,z)&\mapsto E(t,z) =\exp_z(te_1)\end{split}\end{equation*}
is a diffeomorphism.
 Because locally $\phi$ is a function of the sign distance to $\Sigma_{c_0}$, we have for a certain function $f$:
 $$\phi\circ E(t,z)=f(t).$$
 By $\ref{phi1}$, we have
 $$\frac{f'}{f}=ku\, ,$$
 hence $u$ is also a function of $t$.
 We also have 
 $$\Delta t=-\tr \nabla dt=\tr L= (n-p)u+(p-1)\frac{1}{u}.$$
 The equation $\Delta \phi=p(n-1-p)\phi=pk\phi$ implies :
 $$f''-\left((n-p)u+(p-1)\frac{1}{u}\right)f'+kpf=0$$
 We obtain
 $$\frac{f''}{f}=\frac{k+1}{k}\left[\frac{f'}{f}\right]^2-k.$$
 If we let $g(t):=f(t)^{-1/k},$
 we obtain the equation
$$g''-g=0.$$
There are two constants $A, B$ such that
$g(t)=A e^t+Be^{-t}.$
We remark that $AB\not=0$ because $u\not=1$ hence we can find a constant $C$ and a real $\tau$ such that
$$g(t)=\begin{cases}C\cosh(t+\tau) &{\ \rm if}\ AB>0\\
C\sinh(t+\tau) &{\ \rm if}\ AB<0\\
\end{cases}$$

so that 
$$u=\begin{cases}-\tanh(t+\tau) &{\ \rm if}\ AB>0\\
-1/\tanh(t+\tau) &{\ \rm if}\ AB<0\\
\end{cases}$$
 Moreover, because $u$ is always positive we have
 $\tau<0$.
We endow $(-\infty,-\tau)\times S\times T$ with the hyperbolic metric
$$(dt)^2+ \left[f'(t)/f'(0)\right]^2g_{S}+ c_0^{-2}f(t)^2g_{T}$$
so that the map $E$ is a isometry from $(-\delta,\delta)\times S\times T$ onto $\mcO$, the natural extension of this map
 $E(t,z)=\exp_z(te_1)$ becomes an isometric immersion.
The unique continuation property (applied to $\left.\phi\right|_E$) implies again that on $ (-\infty,-\tau)\times S\times T$, 
$$\phi\circ E(t,z)=f(t)^{-k}.$$
As $\phi\circ E$ remains bounded as 
$t$ tends to $-\tau$, hence we  must have $AB>0$.
And the above hyperbolic metric on $(-\infty,-\tau)\times S\times T$ is 
$$(dt)^2+ \left[\frac{\sinh(t+\tau)}{\sinh(\tau)}\right]^2g_{S}+\left[ \frac{\cosh(t+\tau)}{\cosh(\tau)}\right]^2 g_{T}.$$
The metric $\sinh^{-2}(\tau)g_S$ has constant curvature $1$ and the metric $\cosh^{-2}(\tau)g_T$ has constant curvature $-1$, and
\begin{equation}\label{Estarxi}
E^*\xi =C\frac{(\sinh t)^{p+1}}{(\cosh t)^{n-p+1}}dt\wedge d\sigma ;\end{equation}
$d\sigma$ being the volume form of $(S,\sinh^{-2}(\tau)g_S)$. 

But in Fermi coordinate $(0,+\infty)\times \bS^{p-1}\times \H^{n-p}$ around 
a totally geodesic copy of $\H^{n-p}\subset \H^n$, the hyperbolic metric is 
$$(ds)^2+ \sinh^2(s)g_{\bS^{p-1}}+\cosh^2(s) g_{\H^{n-p}}.$$
If $\mcO$ is small enough, we can find an isometry
$\iota \,:\, (-\infty,-\tau)\times S\times T \rightarrow \H^n\setminus \H^{n-p}$.
Eventually, the isometry $\iota\circ E^{-1}$ a priori defined on $\mcO$ can be extended to an isometry $\gamma$ of $\H^n$.
Using this isometry and \ref{Estarxi}, we find the desired expression of $\xi$ on $\mcO$, then the result follows by the unique
continuation property.
\section{Proof of the main theorem}
\subsection{Margulis's decomposition} (See \cite[Chapter D]{BenedettiP} or \cite[\S 12.6]{Ratcliffe})
 Let $X=\Gamma\backslash \H^n$ be a complete hyperbolic manifold and let $\epsilon$ be a positive number
smaller than the Margulis's constant $\epsilon_n$.
The $\epsilon$-thin part of $X$ is the set $X_{<\epsilon}$, where the injectivity radius is smaller that $\epsilon$; we have
$X_{<\epsilon}=V(\Gamma,\epsilon)/\Gamma$ where
$$V(\Gamma,\epsilon):=\{z\in \H^n, \exists \gamma\in \Gamma\setminus\{\mathrm{id}\},\ d(z,\gamma.z)<2\epsilon\}.$$
Let $E$ be the union of all unbounded connected components of $X_{<\epsilon}$ :
$$E=\cup E_j$$
where $\{E_j\}_j$ is the set of the unbounded connected component of $X_{<\epsilon}$.

The topology of such  an $\epsilon$-end is well known. 
When $E_j$ is an unbounded connected component of $X_{<\epsilon}$, there is a point $p_j\in\partial_\infty \H^n$ and a parabolic subgroup
$$\Gamma_j:=\{\gamma\in \Gamma, \gamma.p_j=p_j\},$$
such that on the description of the hyperbolic space as the upper-half space model
 $$\left((0,\infty)\times\R^{n-1}, y^{-2}((dy)^2+(dx)^2) \right)$$
  where the point $p_j$ is at $\infty$, then $\Gamma_j$ acts freely on $\R^{n-1}$ and
  $E_j$ is homeomorphic to $(1,\infty)\times \left(\Gamma_j\backslash \R^{n-1}\right)$.
  Let $F_j$ be the flat manifold $\Gamma_j\backslash \R^{n-1}$ and $S_j\subset F_j$ be a soul of $F_j$, then there is a
   maximal $\Gamma_j$-invariant $r$-plane
  $\tilde S_j\subset \R^{n-1}$ that is the pull-back of a soul by the natural projection
  $\R^{n-1}\rightarrow F_j:= \Gamma_j\backslash \R^{n-1},$ i.e. $S_j=\Gamma_j\backslash \tilde S_j.$
  
  Moreover, there are always positive constants $y_j, r_j$ such that if $\tilde N_j$ is the $r_j$-neighborhood
  of $\tilde S_j\subset \R^{n-1}$ :
  $$\tilde N_j:=\{x\in \R^{n-1}, d(x,\tilde S_j)<r_j\}$$ 
  and $N_j=\Gamma_j\backslash \tilde N_j$, then the inclusion
  $(y_j,\infty)\times N_j\subset E_j$ is a homotopy equivalence.
  We consider $\widehat E_j=[y_j,\infty)\times N_j$ and 
  $\widehat E=\cup_j \widehat{E_j}$.
  We let $$\Sigma_j=\{y_j\}\times N_j\subset \partial \widehat E_j.$$
  And let $\widehat{X}:=\left(X\setminus\widehat{E}\right)\cup\bigcup_j \Sigma_j$, it is a manifold with boundary
  $$\partial \widehat{X}=\bigcup_j \Sigma_j.$$
  We consider the cohomology $$H^\bullet_c(\widehat{X})$$
  of the complex of differential forms $\alpha$ which are smooth on $\widehat{X}$
  and with compact support, that is there is a $R_0$ such that 
  $\mathrm{support}\, \alpha\subset B(o,R_0)$ and 
  $\mathrm{support}\, \alpha\cap\partial \widehat{X}$ is a compact subset of $\partial \widehat{X}=\bigcup_j \Sigma_j$. 
  In particular for
  all but a finite number of $j$ we have $\mathrm{support}\, \alpha\cap\Sigma_j=\emptyset$.
  We have that 
  $$H^\bullet_c(X\setminus X_{<\epsilon})\simeq H^\bullet_c(\widehat{X})$$
  moreover, when $X$ is oriented we have the Poincar\'e duality isomorphism :
  $$H^\bullet(X, X_{<\epsilon})\simeq \left( H^{n-\bullet}_c(X\setminus X_{<\epsilon})\right)^*$$
  The relative cohomology $H^\bullet(X, X_{<\epsilon})$ is also isomorphic to 
  $H^\bullet (\Gamma,\cup_j\Gamma_j, \R)$ 
  which is similar to a cohomology studied by M. Kapovich \cite{kapovich}; more exactly, if 
  $\Pi$ is the union of the $\Gamma_j$'s whose rank is larger or equal to $2$,  he studied
  $H^\bullet(\Gamma,\Pi,\V)$ where 
  $\V$ is a $R\Gamma$-module (for $R$ a commutative ring with unit).
  \subsection{$L^2$ cohomology and harmonic forms}
   We first recall some classical facts on the space of $L^2$ harmonic forms on a complete Riemannian manifold $(X,g)$.

The first one is the Hodge-De Rham-Kodaira orthogonal decomposition :
$$ L^2(\Lambda^{p}T^*X)=\mcH^p(X)\oplus\overline{ dC^\infty _0(\Lambda^{p-1}T^*X)}\oplus \overline{d^*C^\infty _0(\Lambda^{p+1}T^*X)}$$
where the closure are understood for the $L^2$ topology.

The second one is the reduced $L^2$cohomology interpretation of the space of $L^2$ harmonic forms. Let $Z^p_{L^2}(X)$ be the space of weakly closed $L^2$ $p-$forms :
$$Z^p_{L^2}(X):=\left\{\alpha\in L^2(\Lambda^{p}T^*X), d\alpha=0\right\}\ ;$$
By definition,we have
$$Z^p_{L^2}(X)=\left[d^*C^\infty _0(\Lambda^{p+1}T^*X)\right]^\perp=\mcH^p(X)\oplus\overline{ dC^\infty _0(\Lambda^{p-1}T^*X)}.$$
Hence if we introduce the reduced $L^2$-cohomology space : 
$$\H_{L^2}^p(X)\simeq Z^p_{L^2}(X)/\overline{ dC^\infty _0(\Lambda^{p-1}T^*X)}.$$
$$\mcH^p(X)\simeq Z^p_{L^2}(X)/\overline{ dC^\infty _0(\Lambda^{p-1}T^*X)}.$$

We can now describe the natural map from cohomology with compact support to the space of  $L^2$ harmonic forms 
$$H^p_c(X)\rightarrow \mcH^p(X)$$
in two closely related ways.
The first one is induced by the natural inclusions
$$Z^p_c(X):=\left\{\alpha\in C^\infty_c(\Lambda^{p}T^*X), d\alpha=0\right\}\subset Z^p_{L^2}(X) $$
$$ \mathrm{and}\ dC^\infty _0(\Lambda^{p-1}T^*X)\subset \overline{ dC^\infty _0(\Lambda^{p-1}T^*X)}$$
which induces a map
$$H^p_c(X)=\frac{Z^p_c(X)}{ dC^\infty _c(\Lambda^{p-1}T^*X)}\rightarrow \mcH^p(X)\simeq \frac{Z^p_{L^2}(X)}{\overline{ dC^\infty_c(\Lambda^{p-1}T^*X)}}.$$
The second one is induced by the orthogonal projector onto $\mcH^p(X)$ restricted to $Z^p_c(X)$.
This map is zero on $dC^\infty _c(\Lambda^{p-1}T^*X)$ hence induces a map
$H^p_c(X)\rightarrow \mcH^p(X)$.

\subsection{$L^2$ cohomology and "cuspidal" cohomology}
Let  $X=\Gamma\backslash \H^n$ be a complete hyperbolic manifold.
We'll build a natural map from
the cohomology  space $H^p_c(\widehat{X})$ in the space of harmonic $L^2$ $p$-form $\H_{L^2}^p(X)$ for the degrees $p<(n-1)/2$.
The main point is to extend a closed $p-$form with compact support in $\widehat{X}$ to a closed $L^2$-form on $X$.
Let $p_j$ be the projection
$$p_j\,:\, [y_j,\infty)\times N_j\rightarrow \Sigma_j=\{y_j\}\times N_j$$ and
let $\iota_j$ be the inclusion $\Sigma_j\subset \widehat{X}$

Let $\alpha$ be a smooth closed $p$-form with compact support in $\widehat{X}$. We extend $\alpha$ to $X$ by defining 
$\bar\alpha=\alpha$ on $\widehat{X}\subset X$ and
$$\bar\alpha=p_j^*\left(\iota_j^*\alpha\right)\ \mathrm{ on}\ [y_j,\infty)\times N_j. $$ 
It is easy to verify that
$$\|\bar\alpha\|_{L^2([y_j,\infty)\times N_j)}^2=\|\iota_j^*\alpha\|^2_{L^2(\Sigma_j)}\int_{y_j}^\infty y^{-n+2p} dy$$
is finite if $p<(n-1)/2$.
Because $\alpha $ has compact support in $\widehat{X}$, there is only a finite number of $j$ such that 
$\iota_j^*\alpha\not=0$ hence
$$\bar\alpha\in L^2.$$
Remark that by definition, the $\Sigma_j$'s are open, hence $\iota_j^*\alpha$ has compact support in $\Sigma_j$.
Moreover, it is easy to check that $\bar\alpha$ is weakly closed and that
$$\overline{d\alpha}=d\overline{\alpha}\, .$$
Hence we have a well-defined map :
$$ H^p_c(\widehat{X})\rightarrow \H_{L^2}^p(X)\simeq \mcH^p(X)$$
\begin{rem}\label{smooth}
When $c\in H^p_c(\widehat{X})$, we can always find a $\alpha \in c$ such that for each $j$, $\alpha$ has no normal component
on a neighborhood of $\Sigma_j$ : $\nabla d(.,\Sigma_j)\llcorner\alpha=0$, $\alpha$ being moreover closed, this will imply that near
$\Sigma_j$, $\alpha$ is independant of $r=d(.,\Sigma_j)$ (i.e. invariant by the flow of the radial vector field $\nabla d(.,\Sigma_j)$).
 Then the extension $\bar\alpha$ is smooth on $X$.
\end{rem}
\subsection{Remark on the geometrically finite case}\label{GF} When $\Gamma$ is geometrically finite, R. Mazzeo and R. Phillips have computed the $L^2$
cohomology of $X=\Gamma\backslash \H^n$ in terms of the cohomology of a complex of differential forms which vanish
on certain faces of a compactification of $X$ \cite{MP}.

Indeed such a $X$ can be compactified as a manifold$\bar X$ with corner $\partial_r X\cap \partial_c X$ with boundary
$$\partial \bar X=\partial_r X\cup \partial_c X$$
where $\partial_rX$ is the regular boundary of $X$ (the conformally compact boundary of $X$) and
$\partial_cX=\cup_{t=1}^{n-1} \partial_c(t)$ is the cuspidal boundary of $X$, where $\partial_c(t)$ is the union of
the cuspidal face with rank $t$. 
$\widehat{X}$ is homeomorphic to $\bar X\setminus \partial_rX$.
Let $F_p:=\bigcup_{t<n-1-p} \partial_c(t)$. 
When $p<(n-1)/2$, the result of R.Mazzeo and R.Phillips is that
$$\H_{L^2}^p(X)\simeq H^p(\bar X, \partial_rX\cup F_p).$$
We clearly have a map
$$ H^p(\bar X, \partial_rX\cup F_p)\rightarrow H^p(\bar X,\partial_rX)$$
We consider the  long exact sequence :
$$...\rightarrow H^{p-1}(F_p,\partial_r X)\rightarrow H^p(\bar X, \partial_r X\cup F_p)\rightarrow H^p(\bar X,\partial_rX)
\rightarrow H^{p}(F_p,\partial_r X)\rightarrow ...$$
Now if $t=n-1$, then 
$H^k(\partial_c(n-1),\partial_rX)=H^k(\partial_c(n-1))$.
And when $t<n-1$, then
$$H^k(\partial_c(t),\partial_rX)\simeq H^{n-k}(\partial_c(t), \mathrm{o})$$
where $\mathrm{o}$ is the orientation bundle. But the connected component of $\partial_c(t)$ are homotopic to
compact flat manifold of dimension $t$, hence
we have
$$n-k>t\Rightarrow H^k(\partial_c(t),\partial_rX)=\{0\}.$$
Hence the above long exact sequence implies that :
\begin{prop} If $\Gamma$ is a geometrically finite discrete torsion free subgroup of $\mathrm{SO}(n,1)$ and if $p<(n-1)/2$ then
$$\H_{L^2}^p(X)\simeq H^p(\bar X, \partial_r X\cup F_p)\simeq H^p(\bar X,\partial_r X).$$
\end{prop}
Hence we obtain
\begin{cor} If $\Gamma$ is a geometrically finite discrete torsion free subgroup of $\mathrm{SO}(n,1)$   then for $p<(n-1)/2$ the map
 $$H^p_c(\widehat{X})\rightarrow \H_{L^2}^p(X)\simeq \mcH^p(X)$$ is an isomorphism.
\end{cor}

\subsection{The main result}
\begin{thm}
Let $X=\Gamma\backslash \H$ be a hyperbolic manifold and assume that 
$\delta(\Gamma)=n-1-p$ for some integer $p<(n-1)/2$. If the map
$$H^p_c(\widehat{X})\rightarrow \mcH^p(X)$$ is not zero, then
then $\Gamma$ stabilizes  positively a totally geodesic $n-p$-hyperbolic subspace $\H^{n-p}\subset \H^n$ with 
$\vol (\Gamma\backslash \H^{n-p})<\infty$.
\end{thm}
\subsection{Proof of the main result}
\subsubsection{Preliminary}
We assume that $\Gamma\subset \mathrm{SO}(n,1)$ is a discrete torsion free subgroup and that
$\delta(\Gamma)=n-1-p$ with $p$ an integer such that $2p<(n-1)$ and we assume moreover that
 we can find\begin{itemize}
\item  a non zero $L^2$ harmonic $p$-form $\xi$,
\item a closed $p$-form with compact support in $\widehat{X}$, $\alpha$
\end{itemize}
 such that $\xi$ and  $\bar\alpha$ define the same reduced
$L^2$-cohomology class in $\H_{L^2}^p(X).$
That is to say there is a sequence of
smooth $(p-1)$-forms with compact support $(\beta_k)_k$ such that 
$$\xi-\bar\alpha=L^2\!\!\!-\!\!\!\lim_{k\to\infty} d\beta_k.$$

According to \cite[theorem B]{CP}, the spectrum of the Hodge-De Rham Laplacian on the $(p-1)$-forms on $X$ is bounded from below with 
$$\sigma_p:=n-2p.$$
That is we have the spectral gap estimate :
\begin{equation}\label{spectralgap}
\forall \varphi \in C^\infty_c(\Lambda^{p-1}T^*X), \ (n-2p)\|\varphi\|^2_{L^2}\le \|d\varphi\|^2_{L^2}+\|d^*\varphi\|^2_{L^2}=\la\varphi , \Delta \varphi\ra.
\end{equation}

According to the remark (\ref{smooth}) we can always assume that $\bar\alpha$ is smooth. Hence according 
to \cite[prop. 1.3]{nader},
we can find a smooth $(p-1)$-form
$\beta\in L^2(\Lambda^{p-1}T^*X)$ such that 
$$\xi=\bar\alpha+d\beta,\ \mathrm{and }\ d^*\beta=0.$$
Note in particular that this implies $\Delta\beta=(dd^*+d^*d)\beta=d^*\bar\alpha$, so $\Delta\beta$ vanishes
 outside the support of $\bar \alpha$.

\subsubsection{Some estimates}We are going to prove that $\phi:=|\xi|^{\frac{\delta(\Gamma)}{\delta(\Gamma)+1}}$ 
is square integrable.
For this purpose we 'll use Agmon's type estimates as P. Li and J. Wang \cite{LiWang} 
(finite propagation speed argument can also be used) in order to estimate on the growth of $\beta$ and $d\beta$.
There is a finite set $J$ and $R_0>0$ such that 
$$\support\bar \alpha\subset B(o,R_0)\cup\bigcup_{j\in J} p_j^{-1}(N_j)$$
Let $\rho$ be the function distance in $X$ to $B(o,R_0)\cup\bigcup_{j\in J} p_j^{-1}(N_j)$, we have
$$\support\bar\alpha\subset \rho^{-1}\{0\}.$$
For $\tau>0$, we define
$$\rho_\tau=\min(\tau,\rho)\, ;$$ using the fact that $\beta\in L^2$, it is not hard to justify the integration by part formula :
$$ \int_X\left|(d+d^*)\left(e^{\frac{c}{2}\rho_\tau}\beta\right)\right|^2=\int_X \la\Delta\beta,\beta\ra e^{c\rho_\tau}
+ \frac{c^2}{4}\int_{X} |d\rho_\tau|^2\left|\beta\right|^2 e^{c\rho_\tau},$$
remenbering that $\rho$ is zero on the support of $d^* \bar\alpha$, we get :
$$ \int_X\left|(d+d^*)\left(e^{\frac{c}{2}\rho_\tau}\beta\right)\right|^2\le \int_X \la d^* \bar\alpha,\beta\ra 
+\frac{c^2}{4} \int_{\rho^{-1}([0,\tau])} \left|\beta\right|^2 e^{c\rho_\tau}$$
Using the spectral gap estimate (\ref{spectralgap}), we obtain
$$\sigma_p\int_{X} \left|\beta\right|^2 e^{c\rho_\tau} \le \int_X\left|(d+d^*)\left(e^{\frac{c}{2}\rho_\tau}\beta\right)\right|^2,$$
and we easily deduce

$$\left(\sigma_p-\frac{c^2}{4}\right) \int_{\rho^{-1}([0,\tau])} \left|\beta\right|^2 e^{c\rho_\tau}
\le \int_X \la d^* \bar\alpha,\beta\ra .$$
Letting, $\tau$ going to infinity, we obtain the 

\begin{lem}
Let $\sigma_p=n-2p$. Then there is a constant $C$ such that for any $c>2\sqrt{\sigma_p}$ then
\begin{equation}
\label{Agmon}
\left(\sigma_p-\frac{c^2}{4}\right)\int_{X}|\beta|^2(x) e^{c\rho(x)}d\vol(x)\le C.
\end{equation}
\end{lem}
The second estimate of the proof is the following :
\begin{lem}\label{coro1}There is a constant $C$ such that for any $R>0$ :
$$\int_{\rho^{-1}([0,R])}|\beta|^2e^{2\sqrt{\sigma_p\,}\rho(x)} \le C R.$$
\end{lem}

\textit{Proof of the lemma \ref{coro1}.--}
As a matter of fact, according to the inequality (\ref{Agmon}), we have for all
$c\in [0,2\sqrt{\sigma_p})$ :
\begin{equation*}
\begin{split}
\int_{\rho^{-1}([0,R])}|\beta|^2 e^{2\sqrt{\sigma_p\,}\rho}d\vol&\le 
e^{(2\sqrt{\sigma_p\,}-c)R}\int_{\rho^{-1}([0,R])}|\beta|^2 e^{c\rho}d\vol\\
&\le C\left(\sigma_p-\frac{c^2}{4}\right)^{-1} e^{(2\sqrt{\sigma_p\,}-c)R}
\end{split}\end{equation*}
Hence applying this inequality for $c=2\sqrt{\sigma_p}-1/R$ we get
$$\int_{\rho^{-1}([0,R])}|\beta|^2(x) e^{2\sigma_p\rho(x)}d\vol(x)\le C\frac{4R}{4\sqrt{\sigma_p}-1/R}e^{1}$$
This lemma \ref{coro1} implies the following control on the growth of $d\beta$:
\begin{lem}
There is a constant $C$ such that for all $R\ge 1$ :
$$\int_{\rho^{-1}([R,R+1])}|d\beta|^2 \le C Re^{-2\sqrt{\sigma_p} R}.$$
\end{lem}
\proof Let $R\ge 1$,
we use a cut off function $\chi$ such that
$$\support \chi\subset \rho^{-1}[R-1,R+2])\  \mathrm{and}\  \chi=1\  \mathrm{on} \   \rho^{-1}[R,R+1]):$$
and $|d\chi|\le 2$.
We have 
$$\int_{\rho^{-1}([R,R+1])}|d\beta|^2\le \int_X |(d+d^*)(\chi\beta)|^2.$$
Integrating by part, we get :
$$\int_X |(d+d^*)(\chi\beta)|^2=\int_X \chi^2\la\beta,\Delta\beta\ra+\int_X |d\chi|^2|\beta|^2.$$
As $\Delta\beta$ is zero on the support of $\chi$, we get
$$\int_{\rho^{-1}([R,R+1])}|d\beta|^2\le 4\int_{\rho^{-1}([R-1,R+2])}|\beta|^2.$$
\endproof

The last estimate is about the volume growth of the sub-level set of the function $\rho$:
\begin{lem}\label{volume}
There is a constant such that for all $R\ge 1$
$$\vol \left(\rho^{-1}([0,R])\right)\le C e^{(n-1)R} .$$
\end{lem}
\proof We have
 $$ \rho^{-1}([0,R])\subset B(o,R_0+R) \cup \bigcup_{j\in J} \mcV_j(R)$$
 where $\mcV_j(R)$ is the $R$-neighborhood of $p_j^{-1}(N_j)$.
 We can always choose $R_0$ large enough so that for $j\in J$ :
 $\Sigma_j\subset B(o,R_0).$
It is clear that the volume of $B(o,R_0+R)$ satisfies such an estimate.
Now the volume of  $\mcV_j(R)\setminus B(R_0+R)$ is always smaller that the $R+y_j^{-1}R_j$ neighborhood
of $[y_j,\infty)\times S_j$ inside $[y_j,\infty)\times F_j$. Defined
$r:=d(., S_j)$  the distance to the soul $(0,\infty)\times S_j$.
In Fermi coordinate around $(0,\infty)\times S_j$ the Riemannian metric of the 
manifold $(0,\infty)\times F_j$ is 
$$\cosh^2(r)\frac{(dy)^2+(dx)^2}{y^2}+(dr)^2+\sinh^2(r)(d\sigma)^2.$$
Hence the volume of $\mcV_j(R)\setminus B(R_0+R)$ is less than
$$C\int_{y_j}^\infty\frac{dy}{y^n} \int_0^{R+y_j^{-1}R_j} \sinh^{n-1-t_j}(r)\cosh^{t_j}(r) dr\le C e^{(n-1)R},$$
where $t_j=\dim S_j$.
\endproof

\subsubsection{Conclusion}
\begin{lem} if $\delta:=\delta(\Gamma)=n-1-p$,
then the function $\phi:=|\xi|^{\frac{\delta}{\delta+1}}$ is $L^2$.
\end{lem}
\proof 
As a matter of fact : we have
$$\int_{\{\rho\le 1\}}\phi^2\le
 \left(\vol (\{\rho\le 1\})\right)^{\frac{1}{\delta+1}}\left(\int_{\{\rho\le 1\}} |\xi|^2\right)^{\frac{\delta}{\delta+1}}.$$
Moreover for $k\in \N\setminus\{0\}$ :
we have
\begin{equation*}\begin{split}
\int_{\{k\le \rho\le k+1\}}\phi^2&\le
\left(\vol (\{k\le \rho\le k+1\})\right)^{\frac{1}{\delta+1}} \left(\int_{\{k\le \rho\le k+1\}}|\xi|^2\right)^{\frac{\delta}{\delta+1}}\\
 & \le C e^{k\frac{(n-1)}{\delta+1}}k^{\frac{\delta}{\delta+1}}\exp\left(-2\frac{\sqrt{\sigma_p}k\delta}{\delta+1}\right)\\
 &\le C k^{\frac{\delta}{\delta+1}}\exp\left(k\frac{(n-1)-2\sqrt{\sigma_p}\delta}{\delta+1}\right)
 \end{split}\end{equation*}
 But if 
 $p<(n-1)/2$ then
 $$(n-1)-2\sqrt{\sigma_p}\delta=(n-1)-2\sqrt{n-2p}(n-1-p)<0, $$
 hence the result.
\endproof
\begin{rem}\label{compact}
The only place where the hypothesis $p<(n-1)/2$  is used is about the construction of the map 
$$H^p_c(\widehat{X})\rightarrow \H_{L^2}^p(X).$$
However, there is always a natural map from the cohomology of $X$ with compact support in the reduced $L^2$ cohomology.
Our above arguments show that if $\xi\in \mcH^p(X)\setminus\{0\}$ 
with $p<n/2$ has in its $\H_{L^2}^p(X)$ class a representative with compact support, i.e. $\xi$ is the range of the map
$$H^p_c(X)\rightarrow \mcH^p(X)$$ then the function $\phi:=|\xi|^{\frac{\delta}{\delta+1}}$  satisfies :
\begin{enumeroman}
\item If $p<(n-1)/2$, then $\phi\in L^2$
\item If $p=(n-1)/2$, then there is a constant $R$ such that for any $R\ge 1$ :
$$\int_{B(o,R)} \phi^2\le C R^{\frac{2\delta+1}{\delta+1}}.$$
\end{enumeroman}
\end{rem}
We can now finish the proof of the theorem :
According to D. Sullivan's result, the bottom of spectrum of the Laplacian on function on $X$ is $p(n-1-p)$ hence we have the spectral gap
estimate :
$$\forall f\in C^\infty_c(X)\, \ p(n-1-p)\int_X f^2\le \int_X |df|^2.$$

We use a cutoff function
$\chi_R$ such that
$$\mathrm{supp }\chi_R\subset B(o,2R)\ ,\  \chi_R=1\  \mathrm{on}\  B(o,R)\ \mathrm{and\ }
|d\chi_R|\le \frac{2}{R}.$$
Then
\begin{equation*}
\begin{split}
p(n-1-p)\int_X \left|\chi_ R\phi\right|^2&\le\int_X \left|d(\chi_R\phi)\right|^2\\
&\le \int_X \chi_R^2\phi \Delta\phi+\int_X \phi^2 |d\chi_R|^2
\end{split}
\end{equation*}
$\xi$ being closed and co-closed, the Bochner formula and the refined Kato inequality imply that
$$\Delta\phi\le p(n-1-p)\phi$$ hence :
$$\int_{B(o,R)} \phi\left(p(n-1-p)\phi-\Delta\phi\right)\le\int_X \phi^2 |d\chi_R|^2=O\left(R^{-2}\right).$$
Letting $R\to \infty$
we obtain 
$p(n-1-p)\phi-\Delta\phi=0$ everywhere and we have equality everywhere if the refined Kato inequality :
$$\left| \nabla \xi\right|^2=\frac{n+1-p}{n-p}\left|d|\xi|\right|^2$$
We apply our theorem (\ref{thm:Kato}) to $\bar \xi$ the pull back of $\xi$ on $\H^n$.
We notice that $\Gamma$ must stabilize the level set of $|\bar\xi|$. We have two cases :
\begin{enumeroman}
\item In the first case, we have a fundamental domain for the $\Gamma$ of the type
$\{(y,x)\in (0,+\infty)\times \R^{n-1}, x\in D\}$ where $D$ is a fundamental domain for
 the action of $\Gamma$ on $\R^{n-1}\simeq\{y=1\}$. Then $\xi$ can not be in $L^2$.
\item In the second case, $\Gamma$ must stabilize the level set 
$$|\bar\xi|=\sup |\bar\xi|. $$
That is $\Gamma$ stabilizes a totally geodesic copy of $\H^{n-p}\subset \H^n$ and 
$\xi$ being $L^2$, we have
$$\vol (\Gamma\backslash \H^{n-p})<\infty.$$
Eventually, because $\bar \xi$ is $\Gamma$-invariant, the formula given for $\bar \xi$ in the theorem (\ref{thm:Kato}) implies that
 $\Gamma$ acts trivially on the orientation bundle of the normal bundle of $\H^{n-p}\subset \H^n$.
\end{enumeroman}
\subsection{Final remarks}
The above argument and the remark (\ref{compact}) show that we also obtain
 a rigidity result in the case $p=(n-1)/2$ related the cohomology with compact support :
\begin{thm}
Let $X=\Gamma\backslash \H^n$ be a hyperbolic manifold and assume that $p<n/2$, then
\begin{enumeroman}
\item If the critical exponent of $\Gamma$ satisfies $$\delta(\Gamma)<n-1-p$$ then
$X$ carries no non trivial $L^2$ harmonic $p-$form. 
\item If the critical exponent of $\Gamma$ satisfies $$\delta(\Gamma)=n-1-p$$
and if $$\im \left(H^p_c(X)\rightarrow \mcH^p(X)\right)\not= \{0\}$$ then
$\Gamma$ stabilizes cocompactly and positively a totally geodesic $(n-p)$-hyperbolic subspace $\H^{n-p}\subset \H^n$.
\end{enumeroman}
\end{thm}

Eventually, it is true that we always have a injective map \cite{Anderson} :
$$\im(H^p_c(X)\rightarrow H^p(X))\rightarrow \mcH^p(X)$$ hence a corollary of the above rigidity result is :
\begin{cor}Let $X=\Gamma\backslash \H^n$ be a hyperbolic manifold and assume that $p<n/2$.
If the critical exponent of $\Gamma$ satisfies $$\delta(\Gamma)\le n-1-p,$$ then
$$\im(H^p_c(X)\rightarrow H^p(X))=\{0\}.$$ 
\end{cor}

\end{document}